\documentclass[12pt]{article}                
\usepackage{graphicx}
\usepackage{psfrag}
\usepackage{amsmath, amssymb}
\usepackage{mathptmx}

\newcommand{\Av}{{\rm Av}}

\newcommand{\E}{{\cal E}}
\newcommand{\M}{{\cal M}}

\newcommand{\PP}{{\cal P}}
\newcommand{\LL}{{\cal L}}
\newcommand{\eps}{{\varepsilon}}
\newcommand{\bsigma}{{\bar{\sigma}}}

\newcommand{\e}{\mathbb{E}}
\newcommand{\q}{\mathbb{Q}}
\newcommand{\p}{\mathbb{P}}

\newcommand{\Reals}{\mathbb{R}}
\newcommand{\Natural}{\mathbb{N}}

\newcommand\qed{\hfill\hbox{\rlap{$\sqcap$}$\sqcup$}}

\newtheorem{lemma}{Lemma}
\newtheorem{theorem}{Theorem}

\begin{document}

\title{On the replica symmetric solution of the $K$-sat model.}
\author{Dmitry Panchenko\thanks{Texas A\&M University, email: panchenk@math.tamu.edu. Partially supported by NSF grant.}\\
}
\date{}
\maketitle
\begin{abstract} In this paper we translate Talagrand's solution of the $K$-sat model at high temperature into the language of asymptotic Gibbs measures. Using exact cavity equations in the infinite volume limit allows us to remove many technicalities of the inductions on the system size, which clarifies the main ideas of the proof. This approach also yields a larger region of parameters where the system is in a pure state and, in particular, for small connectivity parameter we prove the replica symmetric formula for the free energy at any temperature.
\end{abstract} 
\vspace{0.5cm}

\section{Introduction}

The replica symmetric solution of the random $K$-sat model at high temperature was first proved by Talagrand in \cite{TKsat}, and later the argument was improved in \cite{SG} and, again, in \cite{SG2}. The main technical tool of the proof is the so called cavity method, but there are several other interesting and non-trivial ideas that play an important role. In this paper, we will translate these ideas into the language of asymptotic Gibbs measures developed by the author in \cite{Pspins}. The main advantage of this approach is that the cavity equations become exact in the infinite volume limit, which allows us to bypass all subtle inductions on the size of the system and  to clarify the essential ideas. Using the exact cavity equations, we will also be able to prove that the system is in a pure state for a larger region of parameters. 

Consider an integer $p\geq 2$ and real numbers $\alpha>0,$ called the connectivity parameter, and $\beta>0$, called the inverse temperature parameter. Consider a random function
\begin{equation}
\theta(\sigma_1,\ldots,\sigma_p)=-\beta\prod_{1\leq i\leq p} \frac{1+J_i \sigma_i}{2}
\label{Deftheta}
\end{equation}
on $\{-1,1\}^p$, where $(J_i)_{1\leq i\leq p}$ are independent random signs, $\p(J_i=\pm 1)=1/2.$ Let $(\theta_k)_{k\geq 1}$ be a sequence of independent copies of the function $\theta$, defined in terms of independent copies of $(J_i)_{1\leq i\leq p}$. Using this sequence, we define a Hamiltonian $H_N(\sigma)$ on  $\Sigma_N = \{-1,1\}^N$ by
\begin{equation} 
-H_N(\sigma)=\sum_{k\leq \pi(\alpha N)}\theta_k(\sigma_{i_{1,k}},
\ldots,\sigma_{i_{p,k}}),
\label{Ham}
\end{equation}
where $\pi(\alpha N)$ is a Poisson random variable with the mean $\alpha N$ and the indices $(i_{j,k})_{j,k\geq 1}$ are independent uniform on $\{1,\ldots,N\}$. This is the Hamiltonian of the random $K$-sat model with $K=p$, and our goal will be to compute the limit of the free energy 
\begin{equation}
F_N = \frac{1}{N}\e \log \sum_{\sigma\in \Sigma_N} \exp \bigl(-H_N(\sigma) \bigr)
\end{equation}
as $N\to\infty$ in some region of parameters $(\alpha,\beta)$. It will be convenient to extend the definition of the function $\theta$ from $\{-1,1\}^p$ to $[-1,1]^p$ as follows. Since the product over $1\leq i\leq p$ in (\ref{Deftheta}) takes only two values $0$ and $1$, we can write
$$
\exp \theta(\sigma_1,\ldots,\sigma_p)=1+(e^{-\beta}-1)\prod_{1\leq i\leq p} \frac{1+J_i \sigma_i}{2}.
$$
At some point, we will be averaging $\exp \theta$ over the coordinates $\sigma_1,\ldots,\sigma_p$ independently of each other, so the resulting average will be of the same form with $\sigma_i$ taking values in $[-1,1].$ It will be our choice to represent this average again as $\exp \theta$ with $\theta$ now defined by
\begin{equation}
\theta(\sigma_1,\ldots,\sigma_p)=\log\Bigl(1+(e^{-\beta}-1)\prod_{1\leq i\leq p} \frac{1+J_i \sigma_i}{2}\Bigr).
\label{DefthetaG}
\end{equation}
Of course, on the set $\{-1,1\}^p$ this definition coincides with (\ref{Deftheta}). Note that this function takes values in the interval $[-\beta,0].$

Let us denote by $\Pr[-1,1]$ the set of probability measures on $[-1,1]$. Given $\zeta\in \Pr[-1,1]$, let $(z_i)_{i\geq 1}$ and $(z_{i,j})_{i,j\geq 1}$ be i.i.d. random variables with the distribution $\zeta$ and let
\begin{align}
\PP(\zeta)
=&\,
\log 2 + 
\e \log \Av \exp \sum_{k\leq \pi(p \alpha)} \theta_{k}(z_{1,k},\ldots,z_{p-1,k},\eps)
\nonumber
\\
& -\,
(p-1)\alpha \e \theta(z_1,\ldots,z_p),
\label{PPure}
\end{align}
where $\pi(\alpha p)$ is a Poisson random variable with the mean $\alpha p$ independent of everything else and $\Av$ denotes the average over $\eps\in\{-1,1\}$. The functional $\PP(\zeta)$ is called the replica symmetric formula in this model. Our first result will hold in the region of parameters 
\begin{equation}
\min(4\beta,1) (p-1)p \alpha  <1.
\label{region}
\end{equation}
In this case, we will show that asymptotically the system is always in a pure state in the sense that will be explained in Section \ref{SecPure} and the following holds.
\begin{theorem}\label{ThFE}  If (\ref{region}) holds then
\begin{equation}
\lim_{N\to\infty} F_N = \inf_{\zeta\in \Pr[-1,1]} \PP(\zeta). 
\label{FE}
\end{equation}
\end{theorem}
Notice that when the connectivity parameter $\alpha$ is small, $(p-1)p\alpha<1$, the formula (\ref{FE}) holds for all temperatures, which is a new feature of our approach. One can say more under the additional assumption that
\begin{equation}
\frac{1}{2}(e^\beta-1)(p-1)p \alpha  <1.
\label{region2}
\end{equation}
In particular, in this case one can show that the asymptotic Gibbs measure, which will be defined in the next section, is unique and, as a result, the infimum in (\ref{FE}) can be replaced by $\PP(\zeta)$, where $\zeta$ can be characterized as a fixed point of a certain map arising from the cavity computations.  For $r\geq 1$, let us consider a (random) function $T_r: [-1,1]^{(p-1)r}\to [-1,1]$ defined by
\begin{equation}
T_r\bigl((\sigma_{j,k})_{j\leq p-1, k\leq r}\bigr)
=
\frac{\Av\, \eps \exp A(\eps)}{\Av \exp A(\eps)},
\label{DefTr}
\end{equation}
where 
\begin{equation}
A(\eps) = \sum_{k\leq r} \theta_k(\sigma_{1,k},\ldots,\sigma_{p-1,k}, \eps).
\label{DefTrA}
\end{equation}
We set $T_0 = 0$ and define a map 
\begin{equation}
T: \Pr[-1,1] \to \Pr [-1,1]
\end{equation}
in terms of the functions $(T_r)$ as follows. Given $\zeta \in\Pr[-1,1]$, if we again let $(z_{j,k})_{j\leq p-1, k\geq 1}$ be i.i.d. random variables with the distribution $\zeta$ then $T(\zeta)$ is defined by
\begin{align}
T(\zeta) &= \LL\Bigl(T_{\pi(\alpha p)}\bigl((z_{j,k})_{j\leq p-1, k\leq \pi(\alpha p)}\bigr)\Bigr) 
\nonumber
\\
 &= \sum_{r\geq 0} \frac{(\alpha p)^r}{r!}  e^{-\alpha p }
 \LL\Bigl(T_{r}\bigl((z_{j,k})_{j\leq p-1, k\leq r}\bigr)\Bigr),
\label{DefT}
\end{align}
where $\LL(X)$ denotes the distribution of $X$. In the second line, we simply wrote the distribution as a mixture over possible values of $\pi(\alpha p)$, since this Poisson random variable is independent of everything else. The following is essentially the main result in Chapter 6 in \cite{SG2}.
\begin{theorem}\label{ThFE2}  
If (\ref{region2}) holds then the map $T$ has a unique fixed point, $T(\zeta)=\zeta$. If both (\ref{region}) and (\ref{region2}) hold then $\lim_{N\to\infty} F_N = \PP(\zeta).$
\end{theorem}
As we already mentioned, the main ideas of the proof we give here will be the same as in \cite{SG2} but, hopefully, more transparent. Of course, there is a trade-off in the sense that, instead of working with approximate cavity computations for systems of finite size and using the induction on $N$, one needs to understand how these cavity computations can be written rigorously in the infinite volume limit, which was the main point of \cite{Pspins}. However, we believe that passing through this asymptotic description makes the whole proof less technical and more conceptual. Moreover, the results in  \cite{Pspins} hold for all parameters, and here we simply specialize the general theory to the high temperature region using methods developed in \cite{TKsat, SG, SG2}. 

In the next section, we will review the definition of asymptotic Gibbs measures and recall the main results from \cite{Pspins}, namely, the exact cavity equations and the formula for the free energy in terms of asymptotic Gibbs measures. In Section \ref{SecPure}, we will prove that, under (\ref{region}), all asymptotic Gibbs measures concentrate on one (random) function (so the system is in a pure state) and in Section \ref{SecThFE} we will deduce Theorem \ref{ThFE} from this fact. Finally, in Section \ref{SecThFE2}, we will prove Theorem \ref{ThFE2} by showing that, under (\ref{region}) and (\ref{region2}), the asymptotic Gibbs measure is unique. Of course, as in \cite{SG2}, the same proof works for diluted $p$-spin models as well but, for simplicity of notations, we will work only with the Hamiltonian (\ref{Ham}) of the $p$-sat model.

\section{Asymptotic Gibbs measures}

In this section we will review the main results in \cite{Pspins} starting with the definition of asymptotic Gibbs measures. The Gibbs measure $G_N$ corresponding to the Hamiltonian (\ref{Ham}) is a (random) probability measure on $\{-1,1\}^N$ defined by 
\begin{equation}
G_N(\sigma) = \frac{1}{Z_N} \exp\bigl(- H_N(\sigma)\bigr)
\end{equation} 
where the normalizing factor $Z_N$ is called the partition function. Let $(\sigma^\ell)_{\ell\geq 1}$ be an i.i.d. sequence of replicas drawn from the Gibbs measure $G_N$ and let $\mu_N$ denote the joint distribution of the array of all spins on all replicas, $(\sigma_i^\ell)_{1\leq i\leq N, \ell\geq 1}$, under the average product Gibbs measure $\e G_N^{\otimes \infty}$. In other words, for any choice of signs $a_i^\ell \in\{-1,1\}$ and any $n\geq 1,$
\begin{equation}
\mu_N\Bigl( \bigl\{\sigma_i^\ell = a_i^\ell \ :\ 1\leq i\leq N, 1\leq \ell \leq n \bigr\} \Bigr)
=
\e G_N^{\otimes n}\Bigl( \bigl\{\sigma_i^\ell = a_i^\ell \ :\ 1\leq i\leq N, 1\leq \ell \leq n \bigr\} \Bigr).
\label{muN}
\end{equation}
Let us extend $\mu_N$ to a distribution on $\{-1,1\}^{\Natural\times\Natural}$ simply by setting $\sigma_i^\ell=0$ for $i\geq N+1.$ Let $\M$ be the sets of all possible limits of $(\mu_N)$ over subsequences with respect to weak convergence of measures on the compact product space $\{-1,1\}^{\Natural\times\Natural}$. We will call these limits the \emph{asymptotic Gibbs measures}. One crucial property that these measures inherit from $\mu_N$ is the invariance under the permutation of both spin and replica indices $i$ and $\ell.$ Invariance under the permutation of the replica indices is obvious, and invariance under the permutation of the spin index holds because the distribution of the Hamiltonian (\ref{Ham}) is invariant under any such permutation. In other words, there is symmetry between coordinates in distribution, which is called \emph{symmetry between sites}. 

Because of these symmetries, all asymptotic Gibbs measures have some special structure. By the Aldous-Hoover representation \cite{Aldous, Hoover2}, for any $\mu\in \M$, there exists a measurable function $\sigma:[0,1]^4\to\{-1,1\}$ such that $\mu$ is the distribution of the array 
\begin{equation}
s_i^\ell=\sigma(w,u_\ell,v_i,x_{i,\ell}),
\label{sigma}
\end{equation}
where random variables $w,(u_\ell), (v_i), (x_{i,\ell})$ are i.i.d. uniform on $[0,1]$. The function $\sigma$ is defined uniquely for a given $\mu\in \M$, up to measure-preserving transformations (Theorem 2.1 in \cite{Kallenberg}), so we can identify the distribution $\mu$ of array $(s_i^\ell)$ with $\sigma$. Since, in our case, $\sigma$ take values in $\{-1,1\}$, the distribution $\mu$ is completely encoded by the function
\begin{equation}
\bar{\sigma}(w,u,v) = \e_x \sigma(w,u,v,x)
\label{fop}
\end{equation}
where $\e_x$ is the expectation in $x$ only. The last coordinate $x_{i,\ell}$ in (\ref{sigma}) is independent for all pairs $(i,\ell)$, and we can think of it as flipping a coin with the expected value $\bsigma(w,u_\ell,v_i)$. In fact, given the function (\ref{fop}), we can always redefine $\sigma$ by
$$
\sigma(w,u_\ell,v_i,x_{i,\ell}) = 2 I\Bigl(x_{i,\ell} \leq \frac{1}{2}\bigl(1+ \bsigma(w,u_\ell,v_i) \bigr)\Bigr) -1.
$$
One can think of the function $\bsigma$ in a more geometric way as a Gibbs measure on the space of functions, as follows. It is well known that asymptotically the joint distribution $\mu\in \M$ of all spins contains the same information as the joint distribution of all so called multi-overlaps
\begin{equation}
R_{\ell_1,\ldots, \ell_n}^N = \frac{1}{N} \sum_{1\leq i\leq N} \sigma_i^{\ell_1}\cdots \sigma_i^{\ell_n}
\label{multioverlapN}
\end{equation}
for all $n\geq 1$ and all $\ell_1,\ldots, \ell_n\geq 1$. This is easy to see by expressing the joint moments of one array in terms of the joint moment of the other. In particular, one can check that the asymptotic distribution of the array (\ref{multioverlapN}) over a subsequence of $\mu_N$ converging to $\mu\in \M$ coincides with the distribution of the array 
\begin{equation}
R_{\ell_1,\ldots, \ell_n}
=
\e_v \,\bar{\sigma}(w,u_{\ell_1},v)\cdots \bar{\sigma}(w,u_{\ell_n},v)
\label{multioverlap}
\end{equation}
for $n\geq 1$ and $\ell_1,\ldots, \ell_n\geq 1$, where $\e_v$ denotes the expectation in the last coordinate $v$ only. The average of replicas over spins in (\ref{multioverlapN}) has been replaced by the average of functions over the last coordinate, and we can think of the sequence $(\bsigma(w,u_{\ell_n},\cdot))_{\ell\geq 1}$ as an i.i.d. sequence of replicas sampled from the (random) probability measure 
\begin{equation}
G_w = du \circ \bigl(u\to \bsigma(w,u,\cdot)\bigr)^{-1}
\label{Gibbsw}
\end{equation}
on the space $L^2([0,1], dv) \cap \{ \|\bsigma\|_\infty \leq 1\}$ with the topology of $L^2([0,1], dv)$. Here, both $du$ and $dv$ denote the Lebesgue measure on $[0,1]$. Thus, thanks to the Aldous-Hoover representation, to every asymptotic Gibbs measure $\mu\in \M$ we can associate a function $\bsigma$ on $[0,1]^3$ or a random measure $G_w$ of the above space of functions. One can find a related interpretation in terms of exchangeable random measures in \cite{Austin2}.

The main idea introduced in \cite{Pspins} was a special regularizing perturbation of the Hamiltonian $H_{N}(\sigma)$ that allows to pass some standard cavity computations for the Gibbs measure $G_N$ to the limit and state them in terms of the asymptotic Gibbs measures $\mu\in\M$. We will refer to \cite{Pspins} for details and only mention that the perturbation mimics adding to the system a random number (of order $\log N$) of cavity coordinates from the beginning. Because of this perturbation, treating a finite number of coordinates as cavity coordinates is ``not felt" by the Gibbs measure, which results in a number of useful properties in the limit. The perturbation is small enough and does not affect the limit of the free energy $F_N$. In the rest of this section, we will describe the cavity equations in terms of the functions $\sigma$ in (\ref{sigma}) and state some of their consequences.

Let us introduce some notation. We will often need to pick various sets of different spin coordinates in the array $(s_i^\ell)$ in (\ref{sigma}), and it is quite inconvenient to enumerate them using one index $i\geq 1$. Instead, we will use multi-indices $(i_1,\ldots, i_n)$ for $n\geq 1$ and $i_1,\ldots, i_n\geq 1$ and consider 
\begin{equation}
s_{i_1,\ldots, i_n} = \sigma(w,u, v_{i_1,\ldots, i_n},x_{i_1,\ldots, i_n}),
\label{s}
\end{equation}
where $(v_{i_1,\ldots, i_n}), (x_{i_1,\ldots,i_n})$ are i.i.d. uniform on $[0,1]$. In addition to (\ref{s}), we will need
\begin{equation}
\hat{s}_{i_1,\ldots, i_n} = \sigma(w,u, \hat{v}_{i_1,\ldots, i_n}, \hat{x}_{i_1,\ldots, i_n}),
\label{hats}
\end{equation}
for some independent copies $\hat{v}$ and $\hat{x}$ of the sequences $v$ and $x$. Let $(\theta_{i_1,\ldots, i_n})$ and $(\hat{\theta}_{i_1,\ldots, i_n})$ be i.i.d. copies of the random function $\theta$.

Take arbitrary integer $n, m, q, r\geq 1$ such that $n\leq m.$ The index $q$ will represent the number of replicas selected, $m$ will be the total number of spin coordinates and $n$ will be the number of cavity coordinates. The parameter $r\geq 1$ will index certain terms in the cavity equations that are allowed because of the stability properties of the Hamiltonian (\ref{Ham}); these terms played an important role in \cite{Pspins} and will appear in the formulation of the mains results from \cite{Pspins}, but will not be used throughout this paper after that. For each replica index $\ell\leq q$ we  consider an arbitrary subset of coordinates 
$C_\ell\subseteq \{1,\ldots, m\}$ and split them into cavity and non-cavity coordinates
\begin{equation}
C_\ell^1 = C_\ell\cap \{1,\ldots, n\},\,\,\,
C_\ell^2=C_\ell\cap \{n+1,\ldots,m\}.
\label{C12}
\end{equation}
The following quantities represent the cavity fields for $i\geq 1$,
\begin{equation}
A_i(\eps)=\sum_{k\leq \pi_i(\alpha p)} \theta_{k,i}(s_{1,i,k},\ldots,s_{p-1,i,k},\eps),
\label{Ai}
\end{equation}
where $\eps\in \{-1,1\}$ and $(\pi_i(\alpha p))_{i\geq 1}$ are i.i.d. Poisson random variables with the mean $\alpha p$. Let $\e'$ denote the expectation in $u$ and the sequences $x$ and $\hat{x}$, and $\Av$ denote the average over $(\eps_i)_{i\geq 1}$ in $\{-1,1\}^\Natural$ with respect to the uniform distribution. Define
\begin{align}
U_\ell &=\, \e' \Bigl(\Av \Bigl(\prod_{i\in C_\ell^1} \eps_i \exp \sum_{i\leq n} A_i(\eps_i) \Bigr)
\prod_{i\in C_\ell^2} s_i \exp\sum_{k\leq r} \hat{\theta}_k(\hat{s}_{1,k},\ldots,\hat{s}_{p,k}) \Bigr),
\nonumber
\\
V &=\, \e' \Bigl(\Av \Bigl(\exp \sum_{i\leq n} A_i(\eps_i) \Bigr)
\exp\sum_{k\leq r} \hat{\theta}_k(\hat{s}_{1,k},\ldots,\hat{s}_{p,k}) \Bigr).
\label{Ul}
\end{align}
The following result proved in Theorem $1$ in \cite{Pspins} expresses some standard cavity computations in terms of the asymptotic Gibbs measures.
\begin{theorem}\label{ThSC}
For any $\mu\in \M$ and the corresponding function $\sigma$ in (\ref{sigma}), 
\begin{equation}
\e \prod_{\ell\leq q} \e' \prod_{i\in C_\ell}s_i
=\e \frac{\prod_{\ell\leq q}U_\ell}{V^q}.
\label{SC}
\end{equation}
\end{theorem}
The left hand side can be written using replicas as $\e \prod_{\ell\leq q} \prod_{i\in C_\ell}s_i^\ell$, so it represent an arbitrary joint moment of spins in the array (\ref{sigma}). The right hand side expresses what happens to this joint moment when we treat the first $n$ spins as cavity coordinates. As in \cite{Pspins}, we will denote by $\M_{inv}$ the set of distributions of exchangeable arrays generated by functions $\sigma:[0,1]^4\to\{-1,1\}$ as in (\ref{sigma}) that satisfy the cavity equations (\ref{SC}) for all possible choices of parameters. Theorem \ref{ThSC} shows that $\M\subseteq \M_{inv},$ which was the key to proving the formula for the free energy in terms of asymptotic Gibbs measures. Let us consider the functional
\begin{align}
\PP(\mu)
=&\
\log 2 + \e \log \e' \Av \exp \sum_{k\leq \pi(p \alpha)} \theta_{k}(s_{1,k},\ldots,s_{p-1,k},\eps)
\nonumber
\\
&-\, (p-1)\alpha\, \e\log \e' \exp \theta(s_1,\ldots,s_p).
\label{CalP}
\end{align}
The next result was proved in Theorem $2$ in \cite{Pspins}.
\begin{theorem}\label{ThFEG} The following holds,
\begin{eqnarray}
\lim_{N\to\infty} F_N = 
\inf_{\mu\in \M} \PP(\mu) = 
\inf_{\mu\in \M_{inv}} \PP(\mu).
\label{FEG}
\end{eqnarray}
\end{theorem}
\textbf{Remark.} This result was stated in \cite{Pspins} for even $p\geq 2$ only, where this condition was used in the proof of the Franz-Leone upper bound \cite{FL}. However, in the case of the $p$-sat model the proof works for all $p$ without any changes at all, as was observed in Theorem 6.5.1 in \cite{SG2}. The condition that $p$ is even is needed in the corresponding result for the diluted $p$-spin model, and that is why it appears in \cite{PT, Pspins}, where both models were treated at the same time.

For some applications, it will be convenient to rewrite (\ref{SC}) in a slightly different form. From now on, we will not be using the terms $\hat{\theta}_k$ in (\ref{Ul}), so we will now set $r=0$. Let us consider some function $f(\sigma_1,\sigma_2)$ on $\{-1,1\}^{m\times q}$ of the arguments
\begin{align}
\sigma_1 &=\, (\sigma_1^\ell,\ldots, \sigma_n^\ell)_{1\leq \ell\leq q} \in \{-1,1\}^{n\times q},
\nonumber
\\
\sigma_2 &=\, (\sigma_{n+1}^\ell,\ldots, \sigma_m^\ell)_{1\leq \ell\leq q}\in \{-1,1\}^{(m-n)\times q}.
\label{DefArgs}
\end{align}
For example, if we consider the function 
\begin{equation}
f(\sigma_1,\sigma_2) = \prod_{\ell\leq q} \prod_{i\in C_\ell}\sigma_i^\ell
= \prod_{\ell\leq q} \Bigl(\prod_{i\in C_\ell^1}\sigma_i^\ell \prod_{i\in C_\ell^2}\sigma_i^\ell \Bigr)
\label{Deff}
\end{equation}
then the left hand side of (\ref{SC}) can be written as $\e f(s_1,s_2)$, where $s_1$ and $s_2$ are the corresponding subarrays of $(s_{i}^\ell)$ in (\ref{sigma}). To rewrite the right hand side, similarly to (\ref{s}), let us consider 
\begin{equation}
s_{i_1,\ldots, i_n}^\ell = \sigma(w,u_\ell, v_{i_1,\ldots, i_n},x_{i_1,\ldots, i_n}^{\ell}),
\label{sell}
\end{equation}
where, as always, all the variables are i.i.d. uniform on $[0,1]$ for different indices and define, for $\eps=(\eps_i^\ell)_{i\leq n, \ell\leq q}\in \{-1,1\}^{n\times q}$,
\begin{equation}
\E(\eps) = \prod_{\ell\leq q} \exp \sum_{i\leq n} A_{i,\ell}(\eps_i^l),
\label{DefE0}
\end{equation}
where
\begin{equation}
A_{i,\ell}(\eps_i^\ell) = \sum_{k\leq \pi_i(\alpha p)} \theta_{k,i}(s_{1,i,k}^\ell,\ldots, s_{p-1,i,k}^\ell, \eps_i^\ell).
\end{equation}
Then, with this notation, the equation (\ref{SC}) can be rewritten as
\begin{equation}
\e f(s_1,s_2)
=
\e \frac{\e' \Av f(\eps,s_2) \E(\eps)}{\e' \Av\, \E(\eps)}.
\label{SCinf}
\end{equation}
Simply, we expressed a product of expectations $\e'$ over replicas $\ell\leq q$ by an expectation of the product, using replicas of the random variables $u$ and $x$ that are being averaged. Since any function $f$ on $\{-1,1\}^{m\times q}$ is a linear combination of monomials of the type (\ref{Deff}), (\ref{SCinf}) holds for any such $f$. From here, it is not difficult to conclude that for any functions $f_1,\ldots, f_k$ on $\{-1,1\}^{m\times q}$ and any continuous function $F: \Reals^k\to\Reals,$
\begin{equation}
\e F \bigl(\e' f_1(s_1,s_2),\ldots, \e' f_k(s_1,s_2)\bigr)
=
\e F\Bigl(\frac{\e' \Av f_1(\eps,s_2) \E(\eps)}{\e' \Av\, \E(\eps)},\ldots, 
\frac{\e' \Av f_k(\eps,s_2) \E(\eps)}{\e' \Av\, \E(\eps)}\Bigr).
\label{SCinFk}
\end{equation} 
It is enough to prove this for functions $F(a_1,\ldots,a_k) = a_1^{n_1}\cdots a_k^{n_k}$ for integer powers $n_1,\ldots,n_k\geq 0$, and this immediately follows from (\ref{SCinf}) by considering $f$ on $q(n_1+\ldots+n_k)$ replicas given by the product of copies of $f_1,\ldots, f_k$ on different replicas, so that each $f_i$ appears $n_i$ times in this product.

\section{Pure state}\label{SecPure}

In this section, we will show that in the region (\ref{region}) the function $\bsigma(w,u,v)$ in (\ref{fop}) corresponding to any $\mu\in \M_{inv}$ essentially does not depend on the coordinate $u$. In other words, for almost all $w$, the Gibbs measure $G_w$ in (\ref{Gibbsw}) is concentrated on one function in $L^2([0,1], dv) \cap \{ \|\bsigma\|_\infty \leq 1\}$. This is expressed by saying that the \emph{system is in a pure state}.
\begin{theorem}\label{ThPure}
Under (\ref{region}), $\bsigma(w,u,v) = \e_u \bsigma(w,u,v)$ for almost all $w,u,v \in [0,1]$, where $\e_u$ denotes the expectation in $u$ only.
\end{theorem}
When the system is in a pure state, we will simply omit the coordinate $u$ and write $\bsigma(w,v)$. In this case, a joint moment of finitely many spins,
$$
\e \prod_{i,\ell} s_i^\ell 
= \e\prod_{i,\ell}\bsigma(w,u_i,v_\ell)
= \e\prod_{i,\ell}\bsigma(w, v_\ell),
$$ 
does not depend on replica indices, which means that we can freely change them, for example, $\e s_1^1 s_1^2 s_2^1 s_2^2 = \e s_1^1 s_1^2 s_2^3 s_2^4.$ As in \cite{SG2}, the strategy of the proof will be to show that we can change one replica index at a time,
\begin{equation}
\e s_1^1 \prod_{(i,\ell)\in C} s_i^\ell = \e s_1^{\ell'} \prod_{(i,\ell)\in C} s_i^\ell,
\label{EqPure}
\end{equation}
where a finite set of indices $C$ does not contain $(1,1)$ and $(1,\ell')$. Using this repeatedly, we can make all replica indices different from each other, showing that any joint moment depends only on how many times each spin index $i$ appears in the product. Of course, this implies that
$$
\e \prod_{i,\ell} s_i^\ell 
= \e\prod_{i,\ell} \e_u\bsigma(w,u,v_\ell),
$$ 
so we could replace the function $\bsigma(w,u,v)$ by $\e_u \bsigma(w,u,v)$ without changing the distribution of the array $(s_i^\ell)$. This would be sufficient for our purposes, since we do not really care how the function $\bsigma$ looks like as long as it generates the array of spins $(s_i^\ell)$ with the same distribution. However, it is not difficult to show that, in this case, the function $\bsigma(w,u,v)$ essentially does not depend on $u$ anyway. Let us explain this first.

\medskip
\noindent
\textbf{Proof of Theorem \ref{ThPure}} \emph{(assuming (\ref{EqPure}))}. If (\ref{EqPure}) holds then $\e s_1^1 s_1^2 s_2^1 s_2^2 = \e s_1^1 s_1^2 s_2^3 s_2^4$. This can also be written in terms of the asymptotic overlaps $R_{\ell,\ell'}$ defined in (\ref{multioverlap}) as $$\e R_{1,2}^2 = \e R_{1,2} R_{3,4}.$$ Since $R_{\ell,\ell'}$ is the scalar product in $(L^2[0,1], dv)$ of replicas $\sigma^\ell$ and $\sigma^{\ell'}$ drawn from the asymptotic Gibbs measure $G_w$ in (\ref{Gibbsw}), 
$$
0= \e R_{1,2}^2 - \e R_{1,2} R_{3,4} = \e \mbox{\rm Var}_{G_w}(\sigma^1\cdot \sigma^2), 
$$
which implies that for almost all $w$ the overlap is constant almost surely. Obviously, this can happen only if $G_w$ is concentrated on one function (that may depend on $w$) and this finishes the proof.
\qed 

\medskip
\medskip
\noindent
In the rest of the section we will prove (\ref{EqPure}). The main idea of the proof will be almost identical to Section 6.2 in \cite{SG2}, even though there will be no induction on the system size. One novelty will be that the cavity equations (\ref{SC}) for the asymptotic Gibbs measures will allow us to give a different argument for large values of $\beta$, improving the dependence of the pure state region on the parameters. We will begin with this case, since it is slightly simpler. 

Without loss of generality, we can assume that $\ell'=2$ in (\ref{EqPure}). Given $m,q\geq 1$, for $j=1,2$, let us consider functions
$f_j(\sigma_1,\sigma_2)$ on $\{-1,1\}^{m\times q}$ with $\sigma_1$ and $\sigma_2$ as in (\ref{DefArgs}). We will suppose that
\begin{equation}
0<f_2 \,\mbox{ and }\,
|f_1|\leq f_2.
\label{Compfs}
\end{equation}
Let us fix $n\leq m$ and, as before, we will treat the first $n$ coordinates as cavity coordinates. Consider the map
\begin{equation}
T: \{-1,+1\}^{m\times q} \to \{-1,+1\}^{m\times q}
\end{equation}
that switches the coordinates $(\sigma_1^1,\ldots, \sigma_n^1)$ with $(\sigma_1^2,\ldots, \sigma_n^2)$ and leaves other coordinates untouched. The statement of the following lemma does not involve $\beta$, but it will be used when $\beta$ is large enough.
\begin{lemma}\label{SecPureLem1}
If $(p-1)p\alpha<1$ and the function $f_1$ satisfies $f_1\circ T = -f_1$ then
\begin{equation}
\e \Bigl| \frac{\e' f_1(s_1,s_2)}{\e' f_2(s_1,s_2)}\Bigr| = 0.
\label{EqSecPureLem1}
\end{equation}
\end{lemma}
To see that (\ref{EqSecPureLem1}) implies (\ref{EqPure}) with $\ell'=2$, take $n=1$, $f_2=1$ and $f_1 = 0.5(\sigma_1^1 -\sigma_1^2)\prod_{(i,\ell)\in C} \sigma_i^\ell.$

\medskip
\noindent
\textbf{Proof.}
By (\ref{Compfs}), the function $f_2$ on $\{-1,1\}^{m\times q}$ is strictly separated from $0$, so we can use (\ref{SCinFk}) with $k=2$ and $F(a_1,a_2)=a_1/a_2$ to get
\begin{equation}
\e \Bigl| \frac{\e' f_1(s_1,s_2)}{\e' f_2(s_1,s_2)}\Bigr| 
=
\e \Bigl| \frac{\e' \Av f_1(\eps, s_2) \E(\eps)}{\e' \Av f_2(\eps, s_2) \E(\eps)}\Bigr|.
\label{Lem1cavity}
\end{equation}
Recall that $\Av$ is the average over $\eps = (\eps_i^\ell)_{i\leq n,\ell\leq q} \in \{-1,1\}^{n\times q}$ and
\begin{equation}
\E(\eps) = \prod_{\ell\leq q} \exp \sum_{i\leq n} A_{i,\ell}(\eps_i^\ell),\,\mbox{ where }\,
A_{i,\ell}(\eps_i^\ell) = \sum_{k\leq \pi_i(\alpha p)} \theta_{k,i}(s_{1,i,k}^\ell,\ldots, s_{p-1,i,k}^\ell, \eps_i^\ell).
\label{DefE}
\end{equation}
For a moment, let us fix all the random variables $\pi_i(\alpha p)$ and $\theta_{i,k}$ and let $r: = \sum_{i\leq n} \pi_i(\alpha p).$ Observe right away that if $r=0$ then $\E(\eps) = 1$ and
\begin{equation}
\Av f_1(\eps, s_2) \E(\eps) = \Av f_1(\eps, s_2) =0.
\label{caser0}
\end{equation}
This is because the average $\Av$ does not change if we switch the coordinates $(\eps_1^1,\ldots, \eps_n^1)$ with $(\eps_1^2,\ldots, \eps_n^2)$ (in other words, just rename the coordinates) and, by assumption,
$$
\Av f_1(\eps, s_2)  = \Av \bigl(f_1(\eps, s_2)\circ T\bigr) = -\Av f_1(\eps, s_2).
$$
Now, let us denote the set of all triples $(j,i,k)$ that appear as subscripts in (\ref{DefE}) by 
\begin{equation}
J = \bigl\{ (j,i,k) \, :\,  j\leq p-1, i\leq n, k\leq \pi_i(\alpha p)\bigr\}.
\label{setJdef}
\end{equation}
If we denote by $\tilde{s}_1 = (s_e^\ell)_{e\in J, \ell\leq q}$ all the coordinates of the array $s$ that appear in $\E(\eps)$ then, for $r\geq 1$, we can think of the averages on the right hand side of (\ref{Lem1cavity}) as functions of $s_2$ and $\tilde{s}_1$,
\begin{equation}
\tilde{f}_j = \tilde{f}_j(\tilde{s}_1,s_2) := \Av f_j(\eps, s_2) \E(\eps).
\label{tildefj}
\end{equation}
Even though $s_2$ and $\tilde{s}_1$ are random variables, for simplicity of notation, here we think of them also as variables of the functions $\tilde{f}_j$. First of all, since $|f_1|\leq f_2$,
$$
|\tilde{f}_1 | \leq \Av |f_1(\eps, s_2)| \E(\eps) 
\leq  \Av f_2(\eps, s_2) \E(\eps) = |\tilde{f}_2|.
$$
Similarly to $T$, let $\tilde{T}$ now be the map that switches the vectors of spins $(s_e^1)_{e\in J}$ and $(s_e^2)_{e\in J}$ in $\tilde{s}_1$ corresponding to the first and second replica. Let us show that $\tilde{f}_1 \circ \tilde{T} = -\tilde{f}_1.$ First, we write
$$
\tilde{f}_1 \circ  \tilde{T} = \Av \bigl(f_1(\eps, s_2)\, (\E(\eps)\circ \tilde{T}) \bigr).
$$
As above, we will use that the average $\Av$ does not change if we switch the coordinates $(\eps_1^1,\ldots, \eps_n^1)$ with $(\eps_1^2,\ldots, \eps_n^2)$, so
$$
\tilde{f}_1 \circ  \tilde{T} = \Av \bigl( (f_1(\eps, s_2)\circ T) (\E(\eps)\circ \tilde{T} T) \bigr).
$$
By assumption, $f_1\circ T = - f_1$ and it remains to notice that $\E(\eps)\circ \tilde{T} T = \E(\eps),$ because $ \tilde{T} T$ simply switches all the terms $A_{i,1}$ and $A_{i,2}$ in the definition of $\E(\eps)$. We showed that (\ref{Lem1cavity}) can be rewritten as
\begin{equation}
\e \Bigl| \frac{\e' f_1(s_1,s_2)}{\e' f_2(s_1,s_2)}\Bigr| 
=
\e \Bigl| \frac{\e'  \tilde{f}_1(\tilde{s}_1, s_2)}{\e' \tilde{f}_2(\tilde{s}_1, s_2)}\Bigr|, 
\label{Lem1cavity2}
\end{equation}
and, conditionally on $\pi_i(\alpha p)$ and $\theta_{i,k}$, the pair of functions $\tilde{f}_1, \tilde{f}_2$ satisfies the same properties as $f_1, f_2$. The only difference is that now $n$ is replaced by the cardinality of the set $J$ in (\ref{setJdef}), equal to $(p-1)r$. For a fixed $n$, let us denote by $D(n)$ the supremum of the left hand side of (\ref{Lem1cavity}) over $m\geq n$ and all choices of functions $f_1, f_2$ with the required properties. Then, the equation (\ref{Lem1cavity2}) implies (first, integrating the right hand side conditionally on all $\pi_i(\alpha p)$ and $\theta_{i,k}$)
\begin{equation}
D(n)\leq \e D((p-1)r) = \e D\bigl((p-1)\pi(n \alpha p)\bigr),
\label{dtodn}
\end{equation}
where $\pi(n\alpha p): = r = \sum_{i\leq n} \pi_i(\alpha p)$ is a Poisson random variables with the mean $n\alpha p.$ Recall that, by (\ref{caser0}), $\tilde{f}_1=0$ when $r=0$, so we can set $D(0)=0$. Also, the assumption $|f_1|\leq f_2$ gives that $D(n)\leq 1$ and, thus, $D(n)\leq n.$ Then, (\ref{dtodn}) implies
$$
D(n)\leq \e (p-1)\pi(n \alpha p) = (p-1)p\alpha n.
$$
Using (\ref{dtodn}) repeatedly, we get, by induction on $j\geq 1$, that $D(n)\leq \bigl((p-1)p\alpha\bigr)^j n$. By assumption, $(p-1)p\alpha<1$, so letting $j\to\infty$ proves that $D(n)=0$ for all $n$. This finishes the proof.
\qed

\medskip
\medskip
\noindent
For small values of $\beta$, we will give a slightly different argument, following Section 6.2 in \cite{SG2}.
\begin{lemma}\label{SecPureLem2}
In the notation of Lemma \ref{SecPureLem1}, suppose that $n=1$ and 
\begin{equation} 
(p-1)p \alpha \beta \exp\bigl( 2\beta +\alpha p (e^{2\beta} - 1)\bigr) <1.
\label{SecPureLem2Eq}
\end{equation}
If $f_1\circ T = -f_1$ then (\ref{EqSecPureLem1}) still holds.
\end{lemma}
\textbf{Proof.} The first part of the proof proceeds exactly the same way as in Lemma \ref{SecPureLem1}, and we obtain (\ref{Lem1cavity2}) for the functions $\tilde{f}_1, \tilde{f}_2$ defined in (\ref{tildefj}). Since $n=1$, we can rewrite (\ref{DefE}) as
\begin{equation}
\E(\eps) = \prod_{\ell\leq q} \exp A_{\ell}(\eps_1^\ell),\,\mbox{ where }\,
A_{\ell} = \sum_{k\leq \pi_1(\alpha p)} \theta_{k}(s_{1,k}^\ell,\ldots, s_{p-1,k}^\ell, \eps_1^\ell),
\label{DefE2}
\end{equation}
and the set (\ref{setJdef}) now becomes
\begin{equation}
J = \bigl\{ (j,k) \, :\,  j\leq p-1, k\leq \pi_1(\alpha p)\bigr\}.
\label{setJdef2}
\end{equation}
Its cardinality if $(p-1)r$, where $r= \pi_1(\alpha p).$ Even though we showed that $\tilde{f}_1 \circ \tilde{T} = -\tilde{f}_1$, we can not draw any conclusions yet since the map $T$ switches only one spins in the first and second replicas, while $\tilde{T}$ switches $(p-1)r$ spins $(s_e^1)_{e\in J}$ and $(s_e^2)_{e\in J}$ in $\tilde{s}_1$, of course, conditionally on $\pi_1(\alpha p)$ and $\theta_{k}$. We will decompose $\tilde{f}_1$ into the sum $\tilde{f}_1 = \sum_{e\in J} \tilde{f}_{e},$ where each $\tilde{f}_e$ satisfies $\tilde{f}_e \circ \tilde{T}_e = -\tilde{f}_e$ with some map $\tilde{T}_e$ that switches  $s_e^1$ and $s_e^2$ only. We begin by writing
$$
\tilde{f}_1 = \frac{1}{2}\bigl(\tilde{f}_1 - \tilde{f}_1\circ \tilde{T}\bigr)
=\frac{1}{2}\Bigl(\tilde{f}_1 - \tilde{f}_1\circ \prod_{e\in J}\tilde{T}_e\Bigr).
$$
If we order the set $J$ by some linear order $\leq$ then we can expand this into a telescopic sum,
$$
\tilde{f}_1 =\sum_{e\in J} \frac{1}{2}\Bigl(\tilde{f}_1 \circ \prod_{e'<e}\tilde{T}_{e'} 
- \tilde{f}_1\circ \prod_{e'\leq e}\tilde{T}_{e'}\Bigr).
$$
Then we simply define
$$
\tilde{f}_{e}
:=
\frac{1}{2}\Bigl(\tilde{f}_1 \circ \prod_{e'<e}\tilde{T}_{e'} 
- \tilde{f}_1\circ \prod_{e'\leq e}\tilde{T}_{e'}\Bigr)
$$
and notice that $\tilde{f}_e \circ \tilde{T}_e = -\tilde{f}_e$, since $\tilde{T}_e\tilde{T}_e$ is the identity. Equation (\ref{Lem1cavity2}) implies 
\begin{equation}
\e \Bigl| \frac{\e' f_1(s_1,s_2)}{\e' f_2(s_1,s_2)}\Bigr| 
\leq
\e \sum_{e\in J} \Bigl| \frac{\e'  \tilde{f}_e(\tilde{s}_1, s_2)}{\e' \tilde{f}_2(\tilde{s}_1, s_2)}\Bigr|. 
\label{Lem1cavity3}
\end{equation}
We keep the sum inside the expectation because the set $J$ is random. Recalling the definition of $\tilde{f}_j$ in (\ref{tildefj}), we can write (for simplicity of notation, we will write $\E$ instead of $\E(\eps)$ from now on)
$$
\tilde{f}_{e}(\tilde{s}_1, s_2)
=
\frac{1}{2}\Av \Bigl({f}_1(\eps,s_2) \Bigl(\E\circ \prod_{e'<e}\tilde{T}_{e'} 
- \E\circ \prod_{e'\leq e}\tilde{T}_{e'} \Bigr)\Bigr).
$$
All the maps $\tilde{T}_{e}$ switch coordinates only in the first and second replica. This means that if we write $\E$ defined in (\ref{DefE2}) as $\E = \E' \E''$ where 
$$
\E' =\exp (A_{1}+A_2),\,\,
\E'' = \prod_{3\leq l\leq q} \exp A_{\ell}
$$
then
\begin{equation}
\tilde{f}_{e}(\tilde{s}_1, s_2)
=
\frac{1}{2}\Av \Bigl({f}_1(\eps,s_2) \E'' \Bigl(\E'\circ \prod_{e'<e}\tilde{T}_{e'} 
- \E'\circ \prod_{e'\leq e}\tilde{T}_{e'} \Bigr)\Bigr).
\label{fstwoEs}
\end{equation}
If $e=(j,k)$ then the terms in the last difference only differ in the term $\theta_{k}(s_{1,k}^\ell,\ldots, s_{p-1,k}^\ell, \eps_1^\ell).$ Since $\theta_k \in [-\beta, 0]$ and $A_1+A_2\leq 0$, we can use that $|e^x-e^y|\leq |x-y|$ for $x, y\leq 0$ to get that
$$
\Bigl|\E'\circ \prod_{e'<e}\tilde{T}_{e'} 
- \E'\circ \prod_{e'\leq e}\tilde{T}_{e'} \Bigr|
\leq 2\beta.
$$
Therefore, from (\ref{fstwoEs}) we obtain
$$
|\tilde{f}_{e}(\tilde{s}_1, s_2)|
\leq \beta  \Av \bigl( |{f}_1(\eps,s_2)| \E''\bigr)
\leq \beta  \Av \bigl( {f}_2(\eps,s_2) \E''\bigr).
$$
Similarly, using that $A_1+A_2\in [-2\beta \pi_1(\alpha p),0]$ we get that
$$
\tilde{f}_{2}(\tilde{s}_1, s_2) = \Av \bigl( f_2(\eps,s_2) \E \bigr) = \Av\bigl( f_2(\eps,s_2) \E' \E''\bigr)
\geq  \exp (-2\beta \pi_1(\alpha p)) \Av \bigl( f_2(\eps,s_2) \E''\bigr),
$$
and together the last two inequalities yield
\begin{equation}
|\tilde{f}_{e}(\tilde{s}_1, s_2)|
\leq
 \beta \exp (2\beta \pi_1(\alpha p)) \tilde{f}_{2}(\tilde{s}_1, s_2).
\label{almthr}
\end{equation}
Let $D$ be the supremum of the left hand side of (\ref{Lem1cavity3}) over all pairs of functions $f_1,f_2$ such that $|f_1|\leq f_2$ and $f_1\circ T = -f_1$ under switching one coordinate in the first and second replicas. Then conditionally on $\pi_1(\alpha p)$ and the randomness of all $\theta_k$, each pair $\tilde{f}_e, \tilde{f}_2$ of the right hand side of (\ref{Lem1cavity3}) satisfies (\ref{almthr}), and we showed above that $\tilde{f}_e \circ \tilde{T}_e = -\tilde{f}_e$  under switching one coordinate in the first and second replicas. Therefore, (\ref{Lem1cavity3}) implies that
\begin{equation}
D \leq D\, \e \sum_{e\in J}  \beta \exp (2\beta \pi_1(\alpha p)) 
= D \beta (p-1) \e \pi_1(\alpha p) \exp (2\beta \pi_1(\alpha p)). 
\label{almthr2}
\end{equation}
Even though, formally, this computation was carried out in the case when $\pi_1(\alpha p)\geq 1$, it is still valid when $\pi_1(\alpha p)=0$ because of (\ref{caser0}). Finally, since $\pi_1(\alpha p)$ has the Poisson distribution with the mean $\alpha p$,
$$
\e \pi_1(\alpha p) \exp (2\beta \pi_1(\alpha p))
=
\sum_{k\geq 0} k e^{2\beta k}\frac{(\alpha p)^k}{k!} e^{-\alpha p}
= \alpha p \exp\bigl( 2\beta +\alpha p (e^{2\beta} - 1)\bigr).
$$
The condition (\ref{SecPureLem2Eq}) together with (\ref{almthr2}), obviously, implies that $D=0$ and this finishes the proof.
\qed

\medskip
\medskip
\noindent
To finish the proof of Theorem \ref{ThPure}, it remains to show that the region (\ref{region}) is in the union of the two regions in the preceding lemmas.
\begin{lemma}
If (\ref{region}) holds then either $p(p-1)\alpha <1$ or (\ref{SecPureLem2Eq}) holds.
\end{lemma}
\textbf{Proof.}
If $\beta\geq 1/4$ then $p(p-1)\alpha <1.$ Now, suppose that $\beta\leq 1/4$ and $p(p-1)\alpha \beta<1/4$. First of all, we can bound the left hand side of (\ref{SecPureLem2Eq}) by
$$
(p-1)p \alpha \beta \exp\bigl( 2\beta +\alpha p (e^{2\beta} - 1)\bigr)
<
\frac{1}{4}\exp\bigl( 2\beta +\alpha p (e^{2\beta} - 1)\bigr).
$$
Using that $e^{2\beta} - 1 \leq \sqrt{e} 2\beta$ for $\beta\leq 1/4$ and $p\alpha\beta<1/4$, we can bound the right hand side by
$$
\frac{1}{4}\exp\Bigl(\frac{1}{2} +2\sqrt{e}p \alpha \beta \Bigr)
\leq
\frac{1}{4}\exp\Bigl(\frac{1}{2} +\frac{1}{2}\sqrt{e} \Bigr)\approx 0.94 <1,
$$
and this finishes the proof.
\qed

\section{Inside the pure state}\label{SecThFE}

Suppose now that the system is in a pure state and, for each $\mu\in \M_{inv}$, the corresponding function $\bsigma(w,u,v)$ does not depend on the second coordinate, in which case we will write it as $\bsigma(w,v)$. Let us begin by proving Theorem \ref{ThFE}.

\medskip
\noindent
\textbf{Proof of Theorem \ref{ThFE}.}
When the system is in a pure state, we can rewrite the functional $\PP(\mu)$ in (\ref{CalP}) as follows. First of all, since the expectation $\e'$ is now only in the random variables $x$, which are independent for all spin and replica indices, we can write
$$
\e' \exp \theta(s_1,\ldots,s_p)
=1+(e^{-\beta}-1)\prod_{1\leq i\leq p} \frac{1+J_i \bsigma_i}{2}
=\exp \theta(\bsigma_1,\ldots,\bsigma_p),
$$
where $\bsigma_i=\e's_i = \e' \sigma(w,u,v_i,x_{i}) = \e' \bsigma(w,u,v_i) = \bsigma(w,v_i).$ Similarly,
$$
 \e' \Av \exp \sum_{k\leq \pi(p \alpha)} \theta_{k}(s_{1,k},\ldots,s_{p-1,k},\eps)
 =
\Av \exp \sum_{k\leq \pi(p \alpha)} \theta_{k}(\bsigma_{1,k},\ldots,\bsigma_{p-1,k},\eps),
$$
where $\bsigma_{i,k}=\bsigma(w,v_{i,k})$. Therefore, the functional $\PP(\mu)$ in (\ref{CalP}) can be written as
\begin{align}
\PP(\mu)
= &\,
\log 2 + \e \log \Av \exp \sum_{k\leq \pi(p \alpha)} \theta_{k}(\bsigma_{1,k},\ldots,\bsigma_{p-1,k},\eps)
\nonumber
\\
&- (p-1)\alpha\, \e\theta(\bsigma_1,\ldots,\bsigma_p).
\label{PmuAgain}
\end{align}
Replacing the average over $w\in[0,1]$ by the infimum, this is obviously bigger than
$$
\inf_{w\in [0,1]} \e_v\Bigl(
\log 2 + \log \Av \exp \sum_{k\leq \pi(p \alpha)} \theta_{k}(\eps, \bsigma_{1,k},\ldots,\bsigma_{p-1,k})
- (p-1)\alpha \theta(\bsigma_1,\ldots,\bsigma_p)
\Bigr),
$$
where $\e_v$ is the expectation only in the random variables $(v_i)$ and $(v_{i,k})$. For a fixed $w$, the random variables  $\bsigma_i$ and $\bsigma_{i,k}$ are i.i.d. and, comparing with (\ref{PPure}), this infimum is bigger than $\inf_{\zeta\in \Pr[-1,1]} \PP(\zeta)$. Since this lower bound holds for all $\mu\in \M_{inv}$, Theorem \ref{ThFEG} then implies that
$$
\lim_{N\to\infty} F_N \geq \inf_{\zeta\in \Pr[-1,1]} \PP(\zeta).
$$
The upper bound follow from the Franz-Leone theorem \cite{FL} by considering functions $\bsigma(w,u,v)$ that depend only on the coordinate $v$ (see Section 2.3 in \cite{Pspins}, and also \cite{PT, SG2}). As we mentioned above, it was observed in Theorem 6.5.1 in \cite{SG2} that the upper bound holds for all $p\geq 2$.
\qed

\medskip
\medskip
\noindent
Let us also write down one consequence of the cavity equations (\ref{SC}) for a system in a pure state. Again, let $\bsigma_i=\bsigma(w,v_i)$ and denote $\bsigma_{j,i,k}=\bsigma(w,v_{j,i,k})$. Let
\begin{equation}
\bsigma_{i}' = \frac{\Av\, \eps \exp A_i(\eps)}{\Av \exp A_i(\eps)},
\label{DefTrI}
\end{equation}
where 
\begin{equation}
A_i(\eps) = \sum_{k\leq \pi_i(\alpha p)} \theta_{k,i}(\bsigma_{1,i,k},\ldots,\bsigma_{p-1,i,k}, \eps).
\label{DefTrI2}
\end{equation}
We will now show that the cavity equations (\ref{SC}) imply the following,
\begin{lemma}\label{LemCP}
If the system is in a pure state, for example in the region (\ref{region}), then
\begin{equation}
\bigl(\bsigma_{i} \bigr)_{i\geq 1} \stackrel{d}{=} \bigl(\bsigma_{i}'\bigr)_{i\geq 1}.
\label{SCPure}
\end{equation}
\end{lemma}
\textbf{Proof.} This can be seen as follows. Take $r=0$ and $n=m$ in (\ref{SC}), so all coordinates will be viewed as cavity coordinates. Since the expectation $\e'$ is now only in the random variables $x$, which are independent for all spin and replica indices, as in the proof of Theorem \ref{ThFE} we can write (slightly abusing notation)
$$
U_\ell = \Av \prod_{i\in C_\ell} \eps_i  \exp \sum_{i\leq n} A_i(\eps_i) \,\,\mbox{ and }\,\,
V = \Av \exp \sum_{i\leq n} A_i(\eps_i),
$$ 
where $A_i(\eps)$ are now given by (\ref{DefTrI2}) instead of (\ref{Ai}), i.e. after averaging the random variables $x$. Therefore, $U_\ell/V =  \prod_{i\in C_\ell} \bsigma_i'$. Since $\e' \prod_{i\in C_\ell} s_i = \prod_{i\in C_\ell} \bsigma_i$, (\ref{SC}) becomes
$$
\e \prod_{\ell\leq q} \prod_{i\in C_\ell} \bsigma_i
=
\e \prod_{\ell\leq q} \prod_{i\in C_\ell} \bsigma_i'.
$$
By choosing $q$ and the sets $C_\ell$ so that each index $i$ appears $n_i$ times gives $\e \prod_{i\leq n} \bsigma_i^{n_i} = \e \prod_{i\leq n} \bsigma_i^{\prime\, n_i}$ and this finishes the proof.
\qed

\section{Proof of Theorem \ref{ThFE2}} \label{SecThFE2}

In this section we will prove Theorem \ref{ThFE2} and we begin with the following key estimate. For a moment, we fix the randomness of $(\theta_k)$ and think of $T_r$ defined in (\ref{DefTr}) as a nonrandom function.
\begin{lemma}\label{LemTLip}
The function $T_r$ defined in (\ref{DefTr}) satisfies
\begin{equation}
\bigl| T_r\bigl((\sigma_{j,k})\bigr) - T_r\bigl((\sigma_{j,k}')\bigr)\bigr|
\leq
\frac{1}{2}(e^\beta -1)\sum_{j, k} |\sigma_{j,k} - \sigma_{j,k}'|.
\label{TLip}
\end{equation}
\end{lemma}
\textbf{Proof.}
Let us compute the derivative of $T_r$ with respect to $\sigma_{1,1}.$ If denote the derivative of 
$$
\theta_1(\sigma_{1,1},\ldots,\sigma_{p-1,1}, \eps) 
=
\log\Bigl(1+(e^{-\beta}-1)\frac{1+J_{p,1}\eps}{2} \prod_{1\leq j\leq p-1} \frac{1+J_{j,1}\sigma_{j,1}}{2} \Bigr)$$ 
with respect to $\sigma_{1,1}$ by
$$
\theta_1' = \exp(-\theta_1) (e^{-\beta}-1)\frac{J_{1,1}}{2} \frac{1+J_{p,1}\eps}{2} \prod_{2\leq j\leq p-1} \frac{1+J_{j,1}\sigma_{j,1}}{2} 
$$
then
$$
\frac{\partial T_r}{\partial \sigma_{1,1}}
=
\frac{\Av\, \eps \theta_1'\exp A(\eps)}{\Av \exp A(\eps)}
-
\frac{\Av\, \eps \exp A(\eps)}{\Av \exp A(\eps)}\,
\frac{\Av\, \theta_1' \exp A(\eps)}{\Av \exp A(\eps)}.
$$
Since $\theta_1\in[-\beta,0]$, we see that $J_{1,1}\theta_1'\in[(1-e^{\beta})/2,0]$ and
$$
\Bigl| \theta_1' -   \frac{\Av\, \theta_1' \exp A(\eps)}{\Av \exp A(\eps)}\Bigr| 
\leq \frac{1}{2}(e^\beta-1),
$$
which implies that $|{\partial T_r}/{\partial \sigma_{1,1}}| \leq (e^\beta-1)/2$. The same, obviously, holds for all partial derivatives and this finishes the proof.
\qed

\medskip
\noindent
\emph{Step 1.}
Let us first show that, under (\ref{region2}), there exists unique fixed point $T(\zeta)=\zeta$. The claim will follow from the Banach fixed point theorem once we show that the map $T$ is a contraction with respect to the Wasserstein metric $W(\p,\q)$ on $\Pr[-1,1]$. This metric is defined by
\begin{equation}
W(\p,\q) = \inf \e|z^1-z^2|,
\end{equation}
where the infimum is taken over all pairs $(z^1,z^2)$ with the distribution in the family $M(\p,\q)$ of measures on $[-1,1]^2$ with marginals $\p$ and $\q$. It is well known that this infimum is achieved on some measure $\mu\in M(\p,\q)$. Let $(z_{j,k}^1,z_{j,k}^2)$ be i.i.d. copies for $j\leq p-1$ and $k\geq 1$ with the distribution $\mu$. By (\ref{TLip}) and Wald's identity,
\begin{align*}
&
\e \bigl| T_{\pi(\alpha p)}\bigl((z_{j,k}^1)\bigr) - T_{\pi(\alpha p)} \bigl((z_{j,k}^2)\bigr)\bigr|
\leq
\frac{1}{2}(e^\beta-1) \e \sum_{j\leq p-1, k\leq \pi(\alpha p)} |z_{j,k}^1 - z_{j,k}^2|
\\
& =  \frac{1}{2}(e^\beta-1) (p-1) p \alpha\, \e  |z_{1,1}^1 - z_{1,1}^2| 
=\frac{1}{2}(e^\beta-1)  (p-1) p \alpha  W(\p,\q).
\end{align*}
On the other hand, by the definition (\ref{DefT}), the pair of random variables on the left hand side,
$$
\Bigl( T_{\pi(\alpha p)}\bigl((z_{j,k}^1)\bigr), T_{\pi(\alpha p)} \bigl((z_{j,k}^2) \bigr) \Bigr),
$$
has the distribution in $M(T(\p), T(\q))$ and, therefore,
$$
W\bigl(T(\p), T(\q)\bigr) \leq \frac{1}{2}(e^\beta-1)(p-1) p \alpha W(\p,\q).
$$
The condition (\ref{region2}) implies that the map $T$ is a contraction with respect to $W$. Since the space $(\Pr[-1,1], W)$ is complete, this proves that $T$ has a unique fixed point $\zeta$.

\medskip
\noindent
\emph{Step 2.}
Now, suppose that both (\ref{region}) and (\ref{region2}) hold. Let $\zeta$ be the unique fixed point $T(\zeta)=\zeta$ and let $\bsigma(w,v,u)$ be the function corresponding to a measure $\mu\in \M_{inv}$ in the statement of Theorem \ref{ThFEG}. By Theorem \ref{ThPure}, we know that $\bsigma$ does not depend on $u$ and, therefore, $\bsigma(w,v)$ satisfies Lemma \ref{LemCP}. Recall that $\bsigma_i = \bsigma(w,v_i)$ and let $(z_i)_{i\geq 1}$ be i.i.d. random variables with the distribution $\zeta$. We will now show that
\begin{equation}
\bigl(\bsigma_{i} \bigr)_{i\geq 1} \stackrel{d}{=} \bigl(z_{i}\bigr)_{i\geq 1},
\label{SCPureS}
\end{equation}
which together with (\ref{PmuAgain}) will imply that $\PP(\mu) = \PP(\zeta)$ for all $\mu\in \M_{inv}$, finishing the proof. (By the way, the fact that $(\bsigma_{i})_{i\geq 1}$ are i.i.d. does not mean that the function $\bsigma(w,u)$ does not depend on $w$; it simply means that the distribution of $(\bsigma_{i})_{i\geq 1}$ is independent of $w$.) To show (\ref{SCPureS}), we will again utilize the Wasserstein metric. For any $n\geq 1$, we will denote by $D(n)$ the Wasserstein distance between the distribution of $(\bsigma_{i})_{i\leq n}$ and the distribution of $(z_{i})_{i\leq n}$ (equal to $\zeta^{\otimes n}$) with respect to the metric $d(x,y) = \sum_{i\leq n}|x_i-y_i|$ on $[-1,1]^n$. For any $r=(r_1,\ldots,r_n) \in\Natural^n$ (we assume now that $0\in\Natural$), let us denote 
$$
p_r = \p\bigl( \pi_1(\alpha p)=r_1,\ldots, \pi_n(\alpha p) = r_n\bigr)= \prod_{i\leq n} \frac{(\alpha p)^{r_i}}{r_i!}e^{-\alpha p}.
$$
Since $\zeta = T(\zeta),$ recalling the definition of $T(\zeta)$ in (\ref{DefT}), we get
\begin{equation}
\zeta^{\otimes n} = T(\zeta)^{\otimes n} = \sum_{r\in \Natural^n} p_r\, \bigotimes_{i\leq n}
 \LL\Bigl(T_{r_i}\bigl((z_{j,k})_{j\leq p-1, k\leq r_i}\bigr)\Bigr),
\label{alm1}
\end{equation}
where the random variables $z_{i,k}$ are i.i.d. and have distribution $\zeta$. Next, similarly to (\ref{DefTr}), let us define
\begin{equation}
T_{i,r_i}\bigl(\sigma_{1,k},\ldots,\sigma_{p-1,k}\bigr) = \frac{\Av\, \eps \exp A_i(\eps)}{\Av \exp A_i(\eps)},
\label{DefTrIAg}
\end{equation}
where 
$$
A_i(\eps) = \sum_{k\leq r_i} \theta_{k,i}(\sigma_{1,k},\ldots,\sigma_{p-1,k}, \eps).
$$
In other words, $T_{i,r_i}$ is defined exactly as $T_{r_i}$, only in terms of independent copies $(\theta_{k,i})$ of $(\theta_k)$. Then, Lemma \ref{LemCP} and (\ref{DefTrI}) imply that 
\begin{equation}
\bigl(\bsigma_i\bigr)_{i\leq n} \stackrel{d}{=}
\sum_{r\in \Natural^n} p_r\,
 \LL\Bigl(  \Bigl(T_{i,r_i}\bigl( (\bsigma_{j,i,k})_{j\leq p-1, k\leq r_i} \bigr) \Bigr)_{i\leq n} \Bigr).
\label{alm2}
\end{equation}
Using the fact that $T_{i,r_i}$ are copies of $T_{r_i}$ defined independently over $i$, we can rewrite (\ref{alm1}) by expressing the product measure as a distribution of a vector with independent coordinates,
\begin{equation}
\zeta^{\otimes n} = \sum_{r\in \Natural^n} p_r\, 
 \LL\Bigl(  \Bigl(T_{i,r_i}\bigl( (z_{j,i,k})_{j\leq p-1, k\leq r_i} \bigr) \Bigr)_{i\leq n} \Bigr),
\label{alm12}
\end{equation}
where the random variables $z_{j,i,k}$ are i.i.d. with the distribution $\zeta$. For a given $r\in \Natural^n$, let us denote by $\p_r$ and $\q_r$ the laws on the right hand side of (\ref{alm2}) and (\ref{alm12}). Since the Wasserstein metric satisfies an obvious inequality for convex combinations of measures
\begin{equation}
W\Bigl(\sum_{r\in \Natural^n} p_r \p_r, \sum_{r\in \Natural^n} p_r \q_r\Bigr)
\leq
\sum_{r\in \Natural^n} p_r W \bigl(\p_r,\q_r \bigr),
\label{Wconvex}
\end{equation}
it remains to estimate the distance between $\p_r$ and $\q_r$. By Lemma \ref{LemTLip},
$$
\sum_{i=1}^{n}
\Bigl|
T_{i,r_i}\bigl( (\bsigma_{j,i,k})_{j\leq p-1, k\leq r_i} \bigr) - T_{i,r_i}\bigl( (z_{j,i,k})_{j\leq p-1, k\leq r_i} \bigr)
\Bigr| 
\leq 
\frac{1}{2}(e^\beta-1)
\sum_{i=1}^{n} \sum_{k=1}^{r_i} \sum_{j=1}^{p-1} \bigl|\bsigma_{j,i,k} - z_{j,i,k}\bigr|.
$$
Choosing the vectors $(\bsigma_{j,i,k})$ and $(z_{j,i,k})$ on the right hand side with the optimal joint distribution that achieves the infimum in the definition of Wasserstein distance and taking expectations proves that
$$
W \bigl(\p_r,\q_r \bigr) \leq \frac{1}{2}(e^\beta-1) D\Bigl((p-1)\sum_{i\leq n} r_i\Bigr).
$$
Plugging this into (\ref{Wconvex}) and using (\ref{alm2}) and (\ref{alm12}) proves that
\begin{equation}
D(n) \leq \frac{1}{2}(e^\beta-1) \sum_{r\in \Natural^n} p_r D\Bigl((p-1)\sum_{i\leq n} r_i\Bigr) 
= \frac{1}{2}(e^\beta-1) \e D\bigl( (p-1) \pi(\alpha p n)\bigr),
\label{Dner}
\end{equation}
where $\pi(\alpha p n)$ is a Poisson random variable with the mean $\alpha p n.$ We start with an obvious bound $D(n)\leq 2n$. Then, by induction on $j$, (\ref{Dner}) implies that
$$
D(n) \leq 2 \Bigl(\frac{1}{2}(e^\beta-1) (p-1) p \alpha\Bigr)^j n
$$
for all $j\geq 1$. Letting $j\to\infty$ proves that $D(n)=0$ for all $n$, since we assumed (\ref{region2}), and this finishes the proof.
\qed


\begin{thebibliography}{99}
\footnotesize

\bibitem{Aldous} Aldous, D.: Representations for partially exchangeable arrays of random variables. J. Multivariate Anal.  {11}, no. 4, 581--598 (1981) 

\bibitem{Austin2} Austin, T.: Exchangeable random measures. Preprint, arXiv:1302.2116 (2013)

\bibitem{FL} Franz, S., Leone, M.: Replica bounds for optimization problems and diluted spin systems. J. Statist. Phys. 111, no. 3-4, 535--564 (2003)

\bibitem{Hoover2} Hoover, D. N.: Row-column exchangeability and a generalized model for probability. Exchangeability in probability and statistics (Rome, 1981), pp. 281--291, North-Holland, Amsterdam-New York (1982) 

\bibitem{Kallenberg} Kallenberg, O.: On the representation theorem for exchangeable arrays. J. Multivariate Anal. 30, no. 1, 137--154 (1989)

\bibitem{PT} Panchenko, D., Talagrand, M.: Bounds for diluted mean-fields spin glass models. Probab. Theory Related Fields 130,  no. 3, 319--336 (2004) 

\bibitem{Pspins} Panchenko, D.: Spin glass models from the point of view of spin distributions. To appear in Ann. of Probab.,  arXiv: 1005.2720 (2010)

\bibitem{TKsat} Talagrand, M.: The high temperature case for the random $K$-sat problem. Probab. Theory Related Fields 119, no. 2, 187--212 (2001)

\bibitem{SG} Talagrand, M.:  Spin Glasses: a Challenge for Mathematicians.  Ergebnisse der Mathematik und ihrer Grenzgebiete. 3. Folge A Series of Modern Surveys in Mathematics, Vol. 43. Springer-Verlag (2003) 

\bibitem{SG2} Talagrand, M.: Mean-Field Models for Spin Glasses. Ergebnisse der Mathematik und ihrer Grenzgebiete. 3. Folge A Series of Modern Surveys in Mathematics, Vol. 54, Springer-Verlag (2011) 


\end{thebibliography}
\end{document}